\documentclass[letterpaper, 10 pt, conference]{ieeeconf}  

\IEEEoverridecommandlockouts                              
\title{\LARGE \bf Decentralized Electric Vehicle Charging Control via a Novel Shrunken Primal Multi-Dual Subgradient (SPMDS) Algorithm}

\author{Xiang Huo$^{\dagger}$ and Mingxi Liu$^{\dagger}$
\thanks{$^{\dagger}$X. Huo and M. Liu are with the Department of Electrical and Computer Engineering at the University of Utah, 50 S Central Campus Drive, Salt Lake City, UT, 84112, USA {\tt{\{xiang.huo, mingxi.liu\}@utah.edu}}.}
}

\usepackage{graphicx}
\usepackage{bm}
\usepackage{mathtools}
\usepackage{amsmath,amsfonts,amssymb}

\usepackage{cite}
\usepackage{algorithm,algpseudocode}
\usepackage{url}

\begin{document}

\maketitle
\thispagestyle{empty}
\pagestyle{empty}

\begin{abstract}
The charging processes of a large number of electric vehicles (EVs) require coordination and control for the alleviation of their impacts on the distribution network and for the provision of various grid services. However, the scalability of existing EV charging control paradigms are limited by either the number of EVs or the distribution network dimension, largely impairing EVs' aggregate service capability and applicability. To overcome the scalability barrier, this paper, motivated by the optimal scheduling problem for the valley-filling service, (1) proposes a novel dimension reduction methodology by grouping EVs (primal decision variables) and establishing voltage (global coupled constraints) updating subsets for each EV group in the distribution network and (2) develops a novel decentralized shrunken primal multi-dual subgradient (SPMDS) optimization algorithm to solve this reduced-dimension problem. The proposed SPMDS-based control framework requires no communication between EVs, reduces over $43\%$ of the computational cost in the primal subgradient update, and reduces up to $68\%$ of the computational cost in the dual subgradient update. The efficiency and efficacy of the proposed algorithm are demonstrated through simulations over a modified IEEE 13-bus test feeder and a modified IEEE 123-bus test feeder.

\end{abstract}

\section{Introduction}

The penetration of electric vehicles (EVs) in the US has reached more than 1.3 million by September 2019 \cite{EEI}, and the EV sales growth rate averages an impressive $25\%$  year-over-year since 2013 \cite{AlternativeFuelsDataCenter}. In the meantime, due to the slow and high-cost infrastructure upgrade, the distribution networks are challenged and impacted by the increasing penetration of EVs in terms of the increased energy loss, voltage sags, harmonics, etc. \cite{lopes2010integration,clement2009impact,bronzini2011coordination,qian2010modeling}.
Therefore, appropriate coordinated EV charging control is needed to alleviate the impacts and postpone the high-cost infrastructure upgrade.

In another perspective, it has been proved that the EV charging flexibility (especially the residential charging at night) can be leveraged for a variety of grid services, including peak shaving, valley-filling, frequency regulation, etc. \cite{masoum2011smart,zhang2014coordinating,liu2013decentralized}. In order to explore the latent service capacity of EVs, a variety of control architectures have been investigated by the power community, including centralized, distributed, and decentralized. Centralized control strategies have been widely researched for optimizing the charging and discharging behaviors of EVs to improve grid operation, e.g., reducing operational cost, avoiding excessive voltage drops, and decreasing energy loss \cite{yang2014improved,richardson2011optimal,esmaili2015multi}. Though being theoretically feasible and efficient, the centralized EV charging controls can hardly be widely deployed as they are not scalable w.r.t. both the network dimension and the number of EVs. Moreover, a centralized control structure requires EVs to communicate with and transmit their private information, e.g., state of charge (SOC), battery capacity, and charger capacity, to an aggregator. This information sharing poses potential risks to the customer privacy.

The poor scalability of the centralized control drives the development of distributed EV charging control. In \cite{yang2016distributed}, the total energy cost was minimized under stochastic renewable loads via distributed EV charging control. Considering the limited local computation and communication capacity, Mohammadi \emph{et al.} in \cite{mohammadi2016fully} developed a distributed consensus+innovations method for EV charging cooperation control. The abovementioned approaches, together with other existing distributed methods, indeed alleviate the scalability issue, however, the inevitable peer-to-peer communication poses cyber-security risks to the entire control framework.


To eliminate the peer-to-peer communication and retain the scalability, decentralized (or hybrid) control are attracting more attentions \cite{braun2016distributed,wang2017hybrid,wen2012decentralized}. More recently, the alternating direction method of multipliers (ADMM) has been largely used in decentralized EV charging control for the minimization of separable cost-related objective functions \cite{peng2014distributed,rivera2016distributed,mercurio2013optimal}. However, in more general cases where the objective functions and constraints are coupled and not separable, e.g., valley-filling and power trajectory tracking, the ADMM-based algorithms cannot work effectively. In the state of the art, few theoretical results have attempted to address this issue. In \cite{koshal2011multiuser}, a regularized primal dual subgradient (RPDS) algorithm was developed by approximating the difference between the original optimal solution  and the regularized counterpart. However, the regularization that was introduced to guarantee convergence led to inevitable convergence errors. As an improvement, Liu \emph{et al.} in \cite{liu2017decentralized} developed a shrunken primal dual subgradient (SPDS) algorithm whose convergence does not rely on the regularization and eliminates convergence errors. Admittedly, RPDS and SPDS are more general and are both scalable w.r.t. the number of EVs, however, like other algorithms (no matter distributed or decentralized), the network dimension presents a hidden scalability issue. Specifically, computations in the algorithm iterations must involve either the network connectivity matrix or the adjacent matrix whose dimensions increase proportionally with the distribution network dimension. Hence, implementing these algorithms in large-scale distribution networks will cause memory overflow and may exceed the computing capacity of the on-board controller. 

Motivated by the above observations, this paper is dedicated to designing a decentralized optimization algorithm that is applicable to generic strongly coupled optimization problems and is scalable w.r.t. both the network dimension and the number of primal decision variables. We will verify the efficiency and efficacy of the proposed algorithm through establishing a decentralized EV charging control framework for the overnight valley-filling where the nonseparable objective function and strongly coupled constraints co-exist.

The contribution of this paper is two-fold: (1) This paper for the first time proposes a dimension reduction approach that extracts and preserves critical network information, partitions primal decision variables into groups to reflect their major network impacts, and partitions the network topology to reduce the complexity; (2) A novel decentralized algorithm, shrunken primal multi-dual subgradient (SPMDS), is proposed in cooperation with the dimension reduction method to achieve the scalability w.r.t. both the number of primal decision variables and the network dimension.



\section{Distribution Network Reconstruction}

\subsection{ Radial Distribution Network Model}

Power flow of a radial distribution network can be described by the DistFlow branch equations which only involve the real power, reactive power, and voltage magnitude \cite{baran1989network}. By omitting the line loss and some higher-order terms, a DistFlow model can be linearized to the LinDistFlow model \cite{baran1989optimal}. It has been shown in \cite{bansal2014plug,gan2014convex,zhang2016scalable} that this linearization has negligible impacts on the model accuracy. Hence, in this paper, we adopt the LinDistFlow model to simplify the power flow description and better illustrate the algorithm design. 

Let $\mathcal{N} = \{0,1,\ldots,n \}$ denote the set of nodes in a radial distribution network, where Node 0 is the slack node that maintains its voltage magnitude at a constant $V_0$, then the LinDistFlow model can be written as \cite{baran1989optimal}
\begin{equation}
\bm{V}(T)=\bm{V}_{0}-2 \bm{R} \bm{P}(T)-2 \bm{X} \bm{Q}(T),
\label{eq:1}
\end{equation}
where $\bm{V}(T) \in \mathbb{R}^{n}$ is the global voltage vector at time $T$ consisting of the squared voltage magnitudes of Nodes 1 to $n$, $\bm{V}_0 \in \mathbb{R}^{n}$ denotes a constant vector with all elements equal to $V_0^2$, and $\bm{P}(T) \in \mathbb{R}^{n}$ and $\bm{Q}(T) \in \mathbb{R}^{n}$ denote the real and reactive power consumption from Node 1 to Node $n$ at time $T$, respectively. In addition, $\bm{R}$ and $\bm{X}$ denote the graphical resistance and reactance matrices defined as
\begin{equation} \label{3}
\begin{aligned}
\bm{R} &\in \mathbb{R}^{n \times n}, \quad R_{ij}=\sum_{(i, j) \in \mathbb{E}_{i} \cap \mathbb{E}_{j}} 
r_{ij}, \\
\bm{X}  &\in \mathbb{R}^{n \times n}, \quad 
X_{ij}=\sum_{(i,j) \in \mathbb{E}_{i} \cap \mathbb{E}_{j}} 
x_{ij},
\end{aligned}
\end{equation}
where $r_{ij}$ and $x_{ij}$ are the line resistance and reactance from Node $i$ to Node $j$, respectively, $\mathbb{E}_{i}$ and $\mathbb{E}_{j}$ are the line sets connecting Node 0 and Node $i$, and Node 0 and Node $j$, respectively \cite{farivar2013equilibrium}. 

\subsection{Dimension Reduction Strategy Considering the Topology of the Radial Distribution Network} \label{section_dimension_reduction}

For general optimal power flow or EV charging control problems, decentralized optimization algorithms normally rely on the iterative primal and dual updates \cite{koshal2011multiuser,liu2017decentralized,yin2009nash}. State-of-the-art methods require the complete network topology information (e.g., connectivity matrix, adjacent matrix, or the graphical resistance matrix $\bm{R}$ as in \eqref{3}) to execute the updates. In the primal and dual updating processes, some interim matrices normally have the dimensions that are tens of times of that of the full network dimension, leading to extra requirements on the on-board memory size and computing power. This severely impairs the scalability of those algorithms. To overcome this, it is critical to reduce the dimension of the distribution network. This reduction can be achieved by abstracting and concentrating the key topology information included in the original $\bm{R}$ and $\bm{X}$ matrices as in \eqref{3}. In what follows, we firstly present the general idea of the dimension reduction and then elaborate on the details. 

EVs' charging impacts on the distribution voltage profile $\bm{V}(T)$ can be quantified by the graphical resistance matrix $\bm{R}$, e.g., each element $R_{ij}$ reflects the impact of the aggregated charging power at Node $j$ on the voltage magnitude of Node $i$. Since all $n$ nodes can have EVs connected, it is rational to consider the full dimension of $\bm{R}$ along the horizontal direction (all columns). However, the vertical direction (rows) of $\bm{R}$ can be reduced as EVs at a specific node (or some nodes as a group) only have major impacts on a sub-vector of the original global voltage vector $\bm{V}(T)$. This drives the following partitioning and grouping strategies for the network nodes and EVs.

\subsubsection{Obtaining the commonly-reduced matrix}
The commonly-reduced matrix, denoted by  $\bm{R}^{all}$, aims to extract the critical elements from the graphical resistance matrix $\bm{R}$. Naturally, $\bm{R}^{all}$ should be defined in a way that attenuates the common impacts of the elements in $\bm{R}$. This leads to each element $R_{ij}^{all}$ in $\bm{R}^{all}$ defined as
\begin{equation}
R_{ij}^{all} \triangleq R_{ij}- \frac{\mathcal{G}(\bm{R})}{n^2},
\label{eq:2}
\end{equation}
and the grand sum is define by $\mathcal{G}(\bm{R}) = \sum_{i=1}^{n}\sum_{j=1}^{n} R_{ij}$. This is an information extraction procedure that elements in $\bm{R}$ which have relatively greater impacts on the nodal voltage magnitudes are preserved.

\subsubsection{Obtaining the peak-preserved matrix}

This is another information extraction procedure, aiming at presenting and preserving the \emph{peak} values in each column of $\bm{R}$. To achieve this, the original graphical resistance matrix $\bm{R}$ is concentrated via shrinking each element by the average value of its column, yielding
\begin{equation}
R_{ij}^{col} \triangleq R_{ij}- \frac{\sum_{i=1}^{n} R_{ij}}{n}.
\label{eq:3}
\end{equation}
Defined in this way, $\bm{R}^{col} \in \mathbb{R}^{n \times n}$, the peak-preserved matrix, shows the peak information of the original graphical resistance matrix $\bm{R}$. An example of $\bm{R}^{col}$ for a modified IEEE 13-bus test feeder used in \cite{liu2017decentralized} is shown in the upper middle of Fig. 1.\begin{figure}[!htb]
    \centering
    \includegraphics[width=0.48\textwidth]{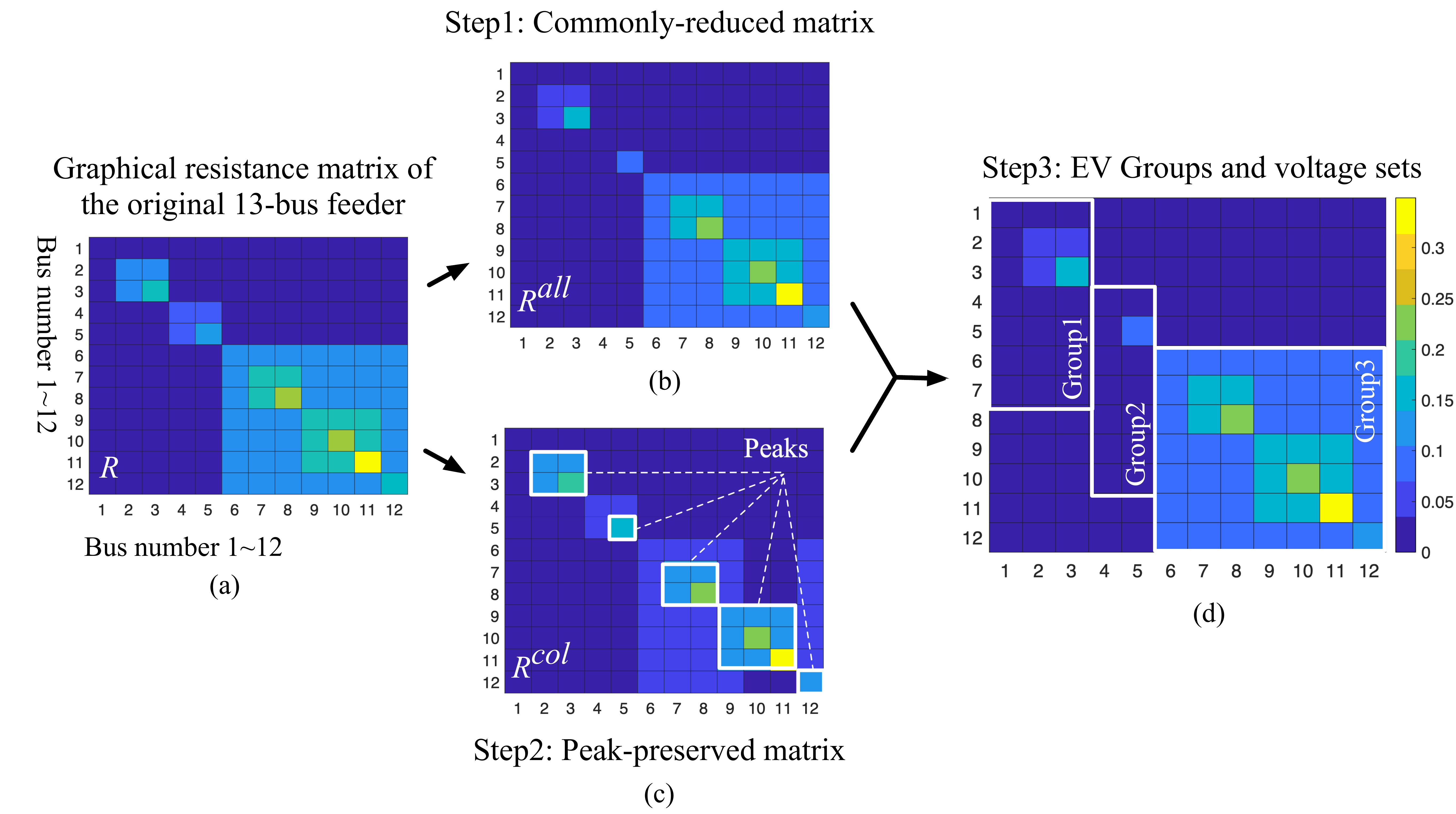}
    \caption{EV grouping and dimension reduction strategy considering the network topology for the modified IEEE 13-bus test feeder.}
    \label{fig:1}
\end{figure} Enabled by the peak-preserved matrix, in total five peaks are preserved from the original distribution network. This procedure ensures that the major impacts on the full voltage magnitude vector $\bm{V}(T)$ by the aggregated charging power at Node $j$ are preserved, i.e., the relatively larger $R_{ij}$ in the $j$th column of $\bm{R}$. Consequently, EVs connected at nodes that have the similar major impacts should be grouped and a voltage sub-vector that reflects the major peak impacts should be constructed correspondingly.


\subsubsection{Dimension reduction}
In this paper, we assume the voltage sub-vector of each EV group has the same dimension. For example, as shown in Fig. 1(d), Group 1 collects EVs connected at Nodes 1-3 and considers voltages at Nodes 1-7, Group 2 collects EVs connected at Nodes 4-5 and considers voltages at Nodes 4-10, and Group 3 collects EVs connected at Nodes 6-12 and considers voltages at Nodes 6-12. The nonuniform voltage partition will be discussed in our future work. Let $d \in \mathbb{Z}_+$ denote the dimension reduction w.r.t. the original global voltage vector $\bm{V}(T)$, then the voltage sub-vector of the $s$th EV group can be written as
\begin{equation} \label{voltage_subvector_def}
    \bm{\hat{V}}_s(T) = \left[\hat{V}_{s,1}(T)~\cdots~\hat{V}_{s,n-d}(T)\right]^{\mathsf{T}}\in \mathbb{R}^{n-d}, 
\end{equation}
for $s=1,\ldots,r$, where $r \in [1,n]$ is the total number of EV groups, and $\hat{V}_{s,l}$, for $l=1,\ldots,n-d$, denotes the $l$th element of the voltage sub-vector of the $s$th EV group. Note that, each voltage sub-vector should include the peaks from $\bm{R}^{col}$ and the preserved elements in $\bm{R}^{all}$ and the voltage sub-vectors of different EV groups may have overlaps.


Let $\tilde{\mathcal{V}}_{1},\ldots,\tilde{\mathcal{V}}_{r}$ denote the voltage subsets and $\tilde{\mathcal{V}}$ denote the full voltage set, then in order to make sure all $n$ nodal voltages are covered by in total $r$ EV groups, the following condition for grouping must be satisfied
\begin{equation}
\tilde{\mathcal{V}}_{1} \cup \tilde{\mathcal{V}}_{2} \cup \cdots \cup \tilde{\mathcal{V}}_{r} = \tilde{\mathcal{V}}.\\
\label{55}
\end{equation}
Further, to ensure \eqref{55} is satisfied, the value of the dimension reduction $d$ must be restricted by
\begin{equation}
0 \leq d \leq n(1-\frac{1}{r}). \\
\label{6}  
\end{equation}

The above three steps extract the key information of the distribution network and reduce the dimension for it. Specifically, Step 1 keeps essential elements in the original matrix; Step 2 preserves the peaks; Step 3 groups EVs and determines the voltage subset for each group based on the commonly-reduced and peak-preserved matrices.

Here, take the modified IEEE 13-bus test feeder for example. Fig. 1 explicitly demonstrates the dimension reduction strategy. The original graphical resistance matrix $\bm{R}$ is presented as a heatmap in Fig. 1(a). Step 1 determines $\bm{R}^{all}$ by extracting the essential elements in $\bm{R}$. The heatmap of $\bm{R}^{all}$ is shown in Fig. 1(b); Step 2 preserves 5 peaks, and the heatmap of $\bm{R}^{col}$ is shown in Fig. 1(c); Step 3 concentrates the essential information in $\bm{R}^{all}$ and $\bm{R}^{col}$ and deploys the dimension reduction rules as shown in Fig. 1(d). Finally, we obtain 3 EV groups and 3 voltage subsets as shown in Fig. 2.
\begin{figure}[!htb]
    \centering
    \includegraphics[width=0.4\textwidth, trim = 0mm 20mm 0mm 0mm, clip]{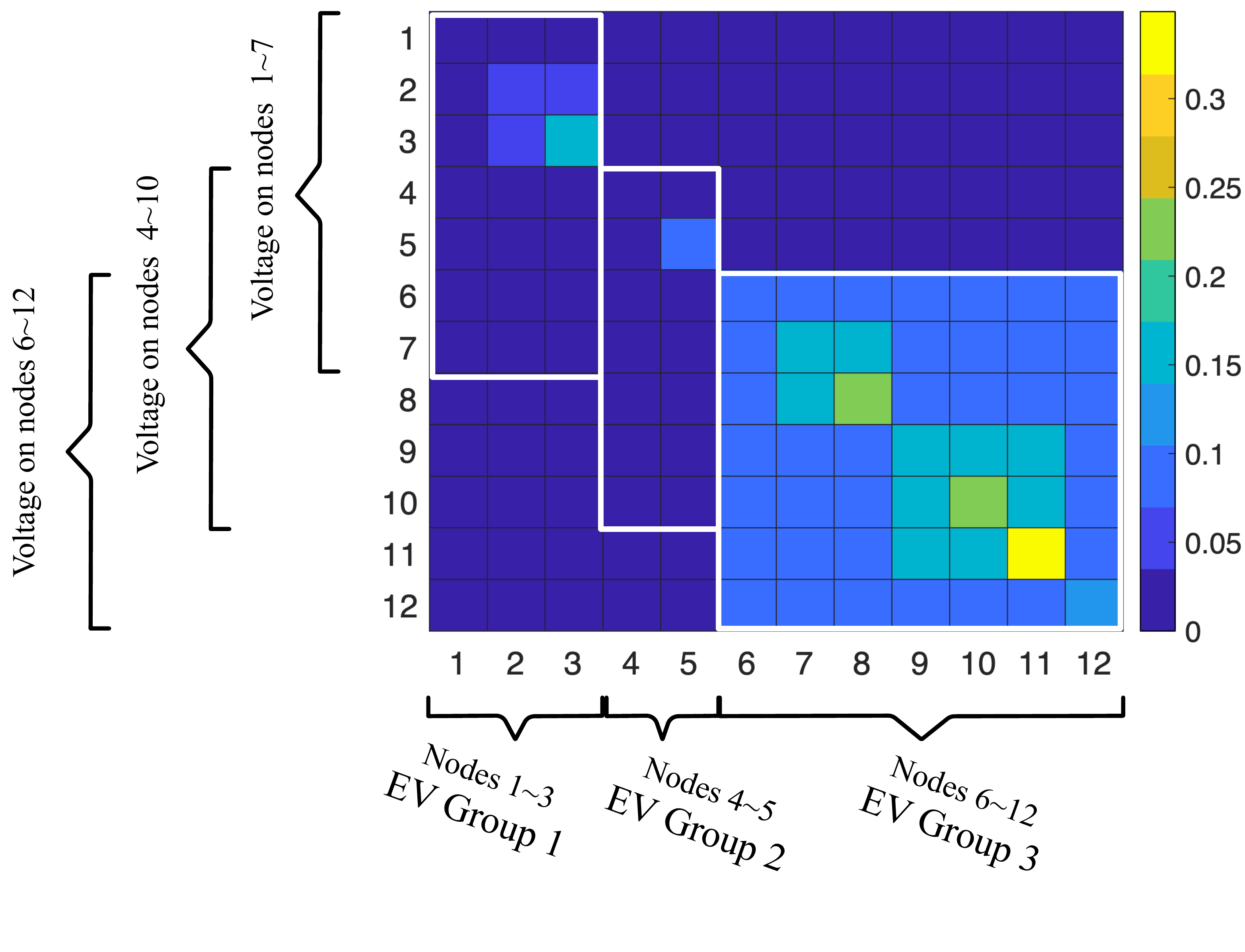}
    \caption{EV grouping and voltage subsets for the modified IEEE 13-bus test feeder (presented via the heatmap of the commonly-reduced matrix).}
    \label{fig:2}
\end{figure}
In this particular case, 500 EVs that are connected to the modified IEEE 13-bus test feeder are divided into three groups, see TABLE \ref{table_example_13} for details.\begin{table}[!htb]
\caption{EV Grouping methodology and voltage subset for each group in the modified IEEE 13-bus test feeder}
\label{table_example_13}
\begin{center}
\begin{tabular}{c|l|l|c}
\hline
Group Name & Location  & Voltage subset & Number of EVs\\
\hline
Group 1 & Nodes 1-3   & Nodes  1-7 & 100\\
Group 2 & Nodes 4-5   & Nodes   4-10& 100\\
Group 3 & Nodes 6-12  & Nodes 6-12 & 300\\
\hline
\end{tabular}
\end{center}
\end{table}
We set the dimension reduction to $d=5$ so as to preserve the complete information in $\bm{R}^{all}$ and $\bm{R}^{col}$. The voltage sub-vectors of Groups 1-3 are
\begin{align*}
 \bm{\hat{V}}_1(T) &= \left[V_{1}(T),\cdots,V_{7}(T)\right]^{\mathsf{T}},\\
 \bm{\hat{V}}_2(T) &= \left[V_{4}(T),\cdots,V_{10}(T)\right]^{\mathsf{T}},\\
 \bm{\hat{V}}_3(T) &= \left[V_{6}(T),\cdots,V_{12}(T)\right]^{\mathsf{T}}.
\end{align*}
Note that $\hat{\bm{V}}_1$, $\hat{\bm{V}}_2$, and $\hat{\bm{V}}_3$ are not unique -- any pair of sub-vectors satisfying the grouping and dimension reduction rules \eqref{eq:2}, \eqref{eq:3}, \eqref{55}, and \eqref{6} can be used. In real charging scenarios, each EV group is responsible for its own voltage subsets, and overlapping between different voltage subsets is allowed.

\section{Decentralized EV Charging Control}
\subsection{Distribution Network and EV Charging Models}
In this paper, we consider that the load at each node is composed of the uncontrollable baseline load and the EV charging load. Without the loss of generality, we assume that the EV charging load is the only controllable load and EVs only consume real power \cite{liu2017decentralized}. Consequently, the LinDistFlow model in \eqref{eq:1} can be rewritten as 
\begin{equation}
\bm{V}(T)=\bm{V}_{0}-\bm{V}_{b}(T)-2 \bm{R} \bm{G}\Bar{\bm{P}}\bm{u}(T),
\label{eq:5}
\end{equation}
where $\bm{V}_{b}(T)$ denotes the voltage drop caused by the baseline load at time $T$, $\bm{u}(T){=}[u_1(T)\cdots u_v(T)]^\mathsf{T} \in \mathbb{R}^{v}$ in $[\bm{0}, \bm{1}]$ contains the normalized charging rates of all $v$ EVs connected at the distribution network. Herein, $\bm{G}$ is the charging aggregation matrix that aggregates the charging power of EVs connected at the same node and $\bm{\Bar{P}}$ is the maximum charging power matrix of all EVs. $\bm{G}$ and $\bm{\Bar{P}}$ are defined by 
\begin{equation}
\bm{G} \triangleq \bigoplus_{j = 1}^{n} \bm{G}_{j} \in \mathbb{R}^{n \times v_j},~\bm{\Bar{P}} \triangleq \bigoplus_{i = 1}^{v} P_{i} \in \mathbb{R}^{v \times v}, \nonumber
\end{equation}
where $\bigoplus$ denotes the matrix direct sum hereinafter, $\bm{G}_j = \bm{1}_{v_j}^{\text{T}}$ is the charging power aggregation vector, ${P}_i$ is the maximum charging power of the $i$th EV, and $v_j$ is the number of EVs connected at node $j$ with $\sum_{j=1}^{n}v_j=v$. 

To simplify the presentation of \eqref{eq:5}, let $\bm{V}_c(T)$ denote $\bm{V}_{0}-\bm{V}_{b}(T)$ and $\bm{D} \in \mathbb{R}^{n \times v}$ denote $2 \bm{R} \bm{G}\Bar{\bm{P}}$, then we have
\begin{equation}
\bm{V}(T)=\bm{V}_c(T) - \bm{D}\bm{u}(T).
\label{eq:6}
\end{equation}

Assume the beginning time of valley-filling is $k$ and the valley-filling process lasts for total $K$ time intervals, then the ending time of valley-filling time is $k+K$. Next, by augmenting the system output described by \eqref{eq:6} along the valley-filling period $[k,k+K]$, we have
\begin{equation}
   \bm{V}_k = \bm{V}_{ck} - \sum_{i=1}^{v} \mathcal{D}_{i} \mathcal{U}_{i,k},
   \label{eq:8}
\end{equation}
where 
\begin{align*}
 &\bm{V}_k = \left[
    \begin{array}{c}
    \bm{V}(k | k) \\ 
    \bm{V}(k+1 | k) \\ 
    \vdots \\ 
    \bm{V}(k+K-1 | k)
    \end{array}
    \right], 
    \bm{V}_{ck}=\left[
    \begin{array}{c}
    \bm{V}_c(k ) \\ 
    \bm{V}_c(k+1) \\ 
    \vdots \\ 
    \bm{V}_c(k+K-1)\end{array}
    \right], \\
& \mathcal{U}_{i,k} {=} \left[\begin{array}{c}
     u_i(k | k) \\ 
     u_i(k+1 | k) \\ 
    \vdots \\ 
    u_i(k+K-1|k)
    \end{array}\right], \mathcal{D}_{i}{=}D_{i} \oplus \cdots \oplus D_{i} {\in} \mathbb{R}^{n K \times K},
\end{align*}
and $D_{i}$ denotes the $i$th column of $\bm{D}$.

Let $x_{i}(T)$ denote the state of charge (SOC) of the $i$th EV at time $T$ and $B_{i}=-\eta_{i} \Delta t \bar{P}_{i}$ denote the maximum charging energy during time $\Delta t$ where $\eta_{i}$ is the charging efficiency, then the charging dynamics of the $i$th EV can be written as 
\begin{equation}
    x_{i}(T+1)=x_{i}(T)+B_{i} u_{i}(T).
\end{equation}
To guarantee all EVs are fully charged by the end of valley-filling, the following equality constraint must be satisfied
\begin{equation}
   \bm{x}(k)+\sum_{i=1}^{v} \mathcal{B}_{i,l} \mathcal{U}_{i}(k) = \bm{0},
   \label{eq:10}
\end{equation}
where  $\bm{x}(k) = \left[x_1(k)~ x_2(k)~\cdots~ x_v(k) \right]^{\mathsf{T}} \in \mathbb{R}^{v}$, $\mathcal{B}_{i, l}=\left[B_{i, c}~ B_{i, c}~ \cdots~ B_{i, c}\right] \in \mathbb{R}^{n \times K}$, and $B_{i, c}$ denotes the $i$th column of the matrix $\bm{B} = \bigoplus_{i = 1}^{v}B_i $.

\subsection{Problem Formulation}
The valley-filling utilizes EVs' charging load to flatten the aggregated total demand, which can be achieved by minimizing the $\ell_2$-norm of the aggregated demand profile \cite{gan2012optimal}. By adopting the same notations as in \cite{liu2017decentralized} and dropping the time indicator $k$ in $\bm{V}_k$, $\bm{V}_{ck}$, and $\mathcal{U}_{i,k}$ hereinafter, we write the objective function of the valley-filling problem as 
\begin{align} 
\mathcal{F}(\mathcal{U}) = \frac{1}{2}\left\|P_{b}+\tilde{P} \mathcal{U}\right\|_{2}^{2},
\label{eq:11}
\end{align}
where $\mathcal{U} =[\mathcal{U}_1^{\mathsf{T}}~ \cdots~ \mathcal{U}_v^{\mathsf{T}}]^{\mathsf{T}}\in \mathbb{R}^{vK}$, $P_{b} \in \mathbb{R}^K$ is the baseline load profile along the valley-filling period, and $\tilde{P} \in \mathbb{R}^{K \times vK}$ is the aggregation matrix for all EVs' charging profiles.

The constraints of the valley-filling problem can be categorized into local constraints and network constraints. For the local constraints, each EV should satisfy
\begin{equation}
  \mathcal{U}_{i} \in \mathbb{U}_{i},
    \mathbb{U}_{i} \triangleq \left\{\mathcal{U}_{i} | \mathbf{0} \leq \mathcal{U}_{i} \leq \mathbf{1}, x_{i}(k)+\mathcal{B}_{i,l} \mathcal{U}_{i}=0\right\}.
    \label{eq:12b}
\end{equation}
For the distribution network, voltage magnitudes at all nodes should be limited within $\left[\underline{v}V_{0}, \bar{v}V_{0}\right]$, where $\underline{v}$ and $\overline{v}$ denote the lower and upper bounds, respectively. With a slight abuse of notation, letting $\bm{V}_0=V_0^2\bm{1}_{nK}$ yields
\begin{equation}
\underline{v}^{2} \bm{V}_{0} \leq \bm{V} \leq \bar{v}^{2} \bm{V}_{0}.\label{eq:13}
\end{equation}
Then by using \eqref{eq:8} and with EVs' charging as the only controllable load, we only need to consider the lower bound for the network voltage constraint, indicating 
\begin{equation}
\bm{V}_c - \sum_{i=1}^{v} \mathcal{D}_{i} \mathcal{U}_{i} \geq \underline{v}^{2} \bm{V}_{0}. \label{eq:14}
\end{equation}
Note that, in more general cases where vehicle-to-grid (V2G) is considered, the upper bound can be added back without affecting the algorithm design. Integrating the objective function \eqref{eq:11} and the constraints \eqref{eq:12b} and \eqref{eq:14}, the valley-filling problem is formulated as
\begin{equation}
\begin{aligned}
& \underset{\mathcal{U}}{\text{min}} & & {\mathcal{F}(\mathcal{U})} \\
& \text{s.t.} & &  \mathcal{U}_{i} \in \mathbb{U}_{i}, \quad \forall i=1,2, \ldots, v,\\
& & &  \underline{v}^{2} \bm{V}_{0} - \bm{V}_c + \sum_{i=1}^{v} \mathcal{D}_{i}\mathcal{U}_{i} \leq \bm{0}. 
\label{eq:15}
\end{aligned}
\end{equation}

\subsection{Shrunken Primal Multi-Dual Subgradient Algorithm}

In this subsection, we develop a scalable approach to solve the optimization problem \eqref{eq:15} in a decentralized fashion. Through Lagrangian relaxation, the problem can be solved by the primal-dual schemes which naturally have the distributed iterate computations across EVs. The relaxed Lagrangian of the problem in \eqref{eq:15} can be written as
\begin{equation}\label{Relaxed_Lag_rep} \mathcal{L}(\mathcal{U}, \bm{\lambda}) =\mathcal{F}(\mathcal{U})+\bm{\lambda}^{\mathsf{T}}\left( \bm{Y}_b + \sum_{i=1}^{v} \mathcal{D}_{i}\mathcal{U}_{i}\right),
\end{equation}
where $\bm{Y}_b = \underline{v}^{2} \bm{V}_{0} - \bm{V}_c$.

The optimization problem in \eqref{eq:15} can be tackled by solving the fixed point problem \cite{boyd2004convex}
\begin{subequations}
\begin{align}
     \mathcal{U}_{i}^{*} &=\Pi_{\mathbb{U}_i}\left(\mathcal{U}_{i}^{*}-\nabla_{\mathcal{U}_{i}} \mathcal{L}\left(\mathcal{U}^{*},\bm{\lambda}^{*}\right)\right),\\
     \bm{\lambda}^{*} &=\Pi_{\mathbb{R}_{+}^{n K}}\left(\bm{\lambda}^{*}+\nabla_{\bm{\lambda}} \mathcal{L}\left(\mathcal{U}^{*}, \bm{\lambda}^{*}\right)\right),
\end{align}
\end{subequations}  
where
\begin{subequations}\label{20}
\begin{align}
    \nabla_{\mathcal{U}_{i}} \mathcal{L}\left(\mathcal{U},\bm{\lambda}\right) &= \tilde{P}^{\mathsf{T}}(P_b+\tilde{P}\mathcal{U}_{i})+\mathcal{D}_{i}^{\mathsf{T}}\bm{\lambda},\label{20a}\\
    \nabla_{\lambda} \mathcal{L}\left(\mathcal{U}, \bm{\lambda}\right) &=\bm{Y}_b + \sum_{i=1}^{v} \mathcal{D}_{i}\mathcal{U}_{i}.\label{20b}
\end{align}
\end{subequations}
Then, by following the dimension reduction method developed in Section \ref{section_dimension_reduction}, we partition EVs into $r$ groups according to the topology of the distribution network. The numbers of EVs in each group are $v_1,\ldots,v_r$, respectively. Use the same $\hat{}$ notation as in \eqref{voltage_subvector_def} to denote the voltage sub-vector, then for the $s$th EV group, the reduced-dimension augmented slack node voltage vector is defined as $\hat{\bm{V}}_{0,s} = V_0^2\bm{1}_{(n-d)K}$ and $\hat{\bm{{V}}}_{c,s} \in \mathbb{R}^{(n-d)K}$ denotes the corresponding sub-vector of $\bm{V}_c$. 
Further, let $\mathcal{D}_{d,i} \in \mathbb{R}^{(n-d) K \times K}$ and $\bm{D}_d \in \mathbb{R}^{(n-d) \times v}$  denote the reduced forms of $\mathcal{D}_{i}$ and $\bm{D}$, respectively (obtained by removing $d$ corresponding rows). Having these definitions, the subgradient calculations in \eqref{20} are modified to
\begin{subequations} \label{modified_gradient}
\begin{align}
    \tilde{\nabla}_{\mathcal{U}_{i}} \mathcal{L}\left(\mathcal{U},\bm{\lambda}\right) 
     &= \tilde{P}^{\mathsf{T}}(P_b+\tilde{P}\mathcal{U}_{i})+\mathcal{D}_{d,i}^{\mathsf{T}}\bm{\lambda}_e, \label{eq:20a}\\
    \tilde{\nabla}_{\bm{\lambda}_s} \mathcal{L}\left(\mathcal{U}, \bm{\lambda}\right)  
     &=\bm{\omega}_s \bm{Y}_{b,s} + \sum_{i=1}^{v} \mathcal{D}_{d,i} \mathcal{U}_{i},\label{eq:20b} \\
    \bm{\lambda}_e &= \bm{\omega}_{1} \odot \bm{\lambda}_1 + \cdots + \bm{\omega}_{r} \odot \bm{\lambda}_r,\label{21c}
\end{align}
\end{subequations}
where $\odot$ denotes elementwise multiplication, $\bm{\lambda}_1,\cdots,\bm{\lambda}_r \in \mathbb{R}^{(n-d)K}$ denote the dual variables of Groups $1$ to $r$, $\bm{\lambda}_e$ denotes the weighted dual variable in the modified primal subgradient \eqref{eq:20a}, $\bm{Y}_{b,s} \in \mathbb{R}^{(n-d)K}$ is a sub-vector of $\bm{Y}_{b}$ in \eqref{Relaxed_Lag_rep}, i.e., 
\begin{equation}
    \bm{Y}_{b,s} = \underline{v}^2\hat{\bm{V}}_{0,s}-\hat{\bm{V}}_{c,s},
\end{equation}
and $\bm{\omega}_1,\bm{\omega}_2,\cdots,\bm{\omega}_r \in \mathbb{R}^{(n-d)K}$ denote the charging impact of each EV group on their designated voltage sub-vector. Take Group 1 in the modified IEEE 13-bus test feeder as an example, the 100 EVs in Group 1 have impacts of \emph{magnitude} $\bm{\omega}_1$ on the voltage sub-vector $\bm{\hat{V}}_1 = \left[V_{1},\ldots,V_{7}\right]^{\mathsf{T}}$, i.e.,
global voltages of the node set \{$1,\ldots,7$\}. Following this definition, $\bm{\omega}_s$, for $s=1,2,\ldots,r$, can be defined as
\begin{equation}\label{eq:22}
\bm{\omega}_s = \left(\sum_{\sum_{j=1}^{s-1}v_j+1}^{\sum_{j=1}^{s}v_j} \mathcal{D}_{d,i} \mathcal{U}_{i} \right) \oslash \left(\sum_{i=1}^{v} \mathcal{D}_{d,i} \mathcal{U}_{i}\right),
\end{equation}
where $\oslash$ denotes the elementwise division. Note that a special case happens when $ \left[\sum_{i=1}^{v} \mathcal{D}_{d,i} \mathcal{U}_{i}\right]_{\jmath} = 0$ for some $\jmath$, where $[\cdot]_{\jmath}$ denotes the $\jmath$th entry. This indicates that no EV is charging at the $\jmath$th time slot. At this point, we define $\left[\bm{\omega}_1\right]_{\jmath} = \left[\bm{\omega}_2\right]_{\jmath} = \ldots = \left[\bm{\omega}_r\right]_{\jmath} = 1$.

Our previous work on SPDS \cite{liu2017decentralized} is dedicated to enabling a decentralized scheme for optimization problems with coupled objective functions and constraints, however its primal and dual updates suffer from large computational load caused by the dimension of the distribution network. In this paper, we propose a novel shrunken primal multi-dual subgradient (SPMDS) algorithm based on the improved primal and dual gradients \eqref{modified_gradient} to circumvent the drawbacks of SPDS. At the $\ell$th iteration, the proposed SPMDS updates the primal and multi-dual variables by following
\begin{subequations}
\begin{align}
\mathcal{U}_{i}^{(\ell+1)} &{=}\Pi_{\mathbb{U}_{i}}\left( \frac{1}{\tau_{\mathcal{U}}} 
\Pi_{\mathbb{U}_{i}}\left( \tau_{\mathcal{U}}\mathcal{U}_{i}^{(\ell)}-\alpha_{i, \ell} \tilde{\nabla}_{\mathcal{U}_{i}} \mathcal{L}\left(\mathcal{U}^{(\ell)}, \bm{\lambda}^{(\ell)}_{e}\right)\right)\right), \label{eq:23a}\\
\bm{\lambda}^{(\ell+1)}_s &{=} \Pi_{\mathbb{D}}\left( \frac{1}{\tau_{\bm{\mathcal{\lambda}}_s}}\Pi_{\mathbb{D}}\left(
\tau_{\bm{\mathcal{\lambda}}_s}
\bm{\lambda}^{(\ell)}_s+\beta_{s,\ell} \tilde{\nabla}_{\bm{\lambda}_s} \mathcal{L}\left(\mathcal{U}^{(\ell)}, \bm{\lambda}^{(\ell)}_s\right)\right)\right), \label{eq:23dual}\\
\bm{\lambda}^{(\ell)}_e &{=} \bm{\omega}_{1}^{(\ell)} \odot \bm{\lambda}^{(\ell)}_1 + \cdots + \bm{\omega}_{r}^{(\ell)} \odot \bm{\lambda}^{(\ell)}_r, \label{eq:23b}
\end{align}
\end{subequations}
for $i=1,\ldots,v$ and $s=1,\ldots,r$, where $\bm{\lambda}^{(\ell)}_1,\bm{\lambda}^{(\ell)}_2, \ldots ,\bm{\lambda}^{(\ell)}_r$ denote the dual variables of Groups $1$ to $r$, respectively, $\mathbb{D}$ is the feasible set of $\bm{\lambda}^{(\ell)}_1,\bm{\lambda}^{(\ell)}_2, \ldots ,\bm{\lambda}^{(\ell)}_r$ (details can be found in \cite{liu2017decentralized}), and $\bm{\lambda}^{(\ell)}_e$ denotes the weighted dual variable of $r$ groups which is commonly adopted by all primal updates \eqref{eq:23a}. In addition, $\tau_{\mathcal{U}}\in [0,1]$ and $\tau_{\bm{\mathcal{\lambda}}_1},\ldots,\tau_{\bm{\mathcal{\lambda}}_r} \in [0,1]$ denote the shrinking parameters of the primal and $r$ dual updates, respectively. $\alpha_{i, \ell} > 0$, for $i=1,2,\ldots,v$, denotes the primal step size of the $i$th EV, and 
$\beta_{s,\ell}>0$, for $s=1,2,\ldots,r$, denotes the dual step size of the $s$th EV group. The convergence criterion is set correspondingly to the convergence of the charging profiles of all EVs, i.e.,
\begin{equation}
    \epsilon=\left\|{\mathcal{U}}^{(\ell+1)}-\mathcal{U}^{(\ell)}\right\|_{2} < \epsilon_0,
    \label{eq:24}
\end{equation}
where $\epsilon_0$ denotes the predefined convergence tolerance. The complete SPMDS algorithm is summarized in Algorithm \ref{SPMDS_alrogithm_1}.
\begin{algorithm} 
\caption{SPMDS Algorithm}
\label{SPMDS_alrogithm_1}
\begin{algorithmic}[1]
\State Calculate the commonly-reduced matrix in \eqref{eq:2}.

\State Calculate the peak-preserved matrix in \eqref{eq:3}.

\State Determine $r$ EV groups and the corresponding voltage subsets $\tilde{\mathcal{V}}_1,\ldots,\tilde{\mathcal{V}}_r$ by \eqref{eq:2}, \eqref{eq:3}, \eqref{55} and \eqref{6}.

\State Parameters initialization: Operator initializes $\bm{\lambda}^{(0)}_{1},\ldots,\bm{\lambda}^{(0)}_{r}$; EVs initialize $\mathcal{U}_i^{(0)}$; Tolerance $\epsilon_0$; Primal and dual step length of $\alpha_{i, \ell}$ and 
$\beta_{1,\ell}>0, \ldots, \beta_{r,\ell}>0$;  Shrinking parameters $\tau_{\mathcal{U}}$ and
$\tau_{\bm{\lambda}_{1}},\ldots,\tau_{\bm{\lambda}_{r}}$; Iteration counter $\ell=0$; Maximum iteration number $\ell_{max}$.

\While{$\epsilon > \epsilon_0$ and $\ell < \ell_{max}$}

\State The central operator computes $\bm{\omega}_{1}^{(\ell)},\ldots,\bm{\omega}_{r}^{(\ell)}$ in \eqref{eq:22}, and calculates and broadcasts $\bm{\lambda}_{e}^{(\ell)}$ in \eqref{eq:23b} to all EVs.

\State Individual EV solves \eqref{eq:23a} and uploads its charging status $\mathcal{U}_i^{(\ell+1)}$ to the operator.

\State The central operator solves \eqref{eq:23dual} to obtain $\bm{\lambda}_{1}^{(\ell+1)},\ldots,\bm{\lambda}_{r}^{(\ell+1)}$.

\State  $\epsilon=\left\|{\mathcal{U}}^{(\ell+1)}-\mathcal{U}^{(\ell)}\right\|_{2}$.

\State $\ell=\ell+1$.

\EndWhile
\end{algorithmic}
\end{algorithm}

\noindent {\bf{Remark 1:}} Though both SPMDS and SPDS emulate the public key encryption, where the dual variables $\bm{\lambda}_e$ and $\bm{\lambda}$ are the public keys while $\mathcal{D}_{d,i}$  and $\mathcal{D}_i$ are the private keys, SPMDS enhances the cyber-security of SPDS. In SPDS, once both $\bm{\lambda}$ and $\mathcal{D}_i$ are sniffed by cyber-attackers, they can be used to reverse engineer the distribution network topology and configuration. In contrast, this reverse engineering cannot be done in SPMDS with the weighted reduced-dimension $\bm{\lambda}_e$ and the reduced-dimension $\mathcal{D}_{d,i}$. \hfill $\blacksquare$

\subsection{Computational Load Reduction Analysis}

SPMDS and SPDS have the similar primal-dual structure, but the former requires less memory and less computational cost for each EV charger (or equivalent onboard controller), and less computation time for the central operator. In what follows, we analyze the computational cost reduction enabled by SPMDS. 

The computation cost difference between \eqref{20a} and  \eqref{eq:20a} lies in the last term $\mathcal{D}_{i}^{\mathsf{T}}\bm{\lambda}$ in \eqref{20a}. Without the EV grouping and dimension reduction strategy, $\mathcal{D}_{i}^{\mathsf{T}}\bm{\lambda}$ requires $(2nK-1)K$ FLOPS from each EV; after the dimension reduction of $d$ for each group, the FLOPS reduce to $2(n-d)K^2-K$. Both \eqref{20a} and \eqref{eq:20a} have the same FLOPS of $4K^2-K$ for the first term $\tilde{P}^{\mathsf{T}}(P_b+\tilde{P}\mathcal{U}_{i})$. Therefore the total FLOPS reduction and the corresponding reduction ratio of the primal gradient  calculation in \eqref{eq:20a} are 
\begin{subequations}
\begin{align}
    \mathcal{F}_{pt} &= (2nK-1)K - ((2n-d)K-1)K  \nonumber\\
    &= 2dK^2,\label{eq:25a}\\
    \mathcal{F}_{pr} &= \frac{2dK^2}{(2nK-1)K+4K^2-K} \nonumber \\
    &= \frac{dK}{(n+2)K-1}.\label{eq:25b}
    \end{align}
\end{subequations}
Eqn. \eqref{eq:25a} shows a proportional relation between the absolute computational cost reduction and the dimension reduction $d$.

For the dual update \eqref{eq:23b}, $\bm{\lambda}_1,\cdots,\bm{\lambda}_r$ can be calculated in parallel. Therefore, the total calculation time for the dual update depends on (i) the time consumed by the EV group that has the largest number of EVs, where we assume this group to be  $r_{m}$ with $v_{m}$ EVs, and (ii) the calculation time for \eqref{eq:20b}. The FLOPS in the dual subgradient calculation in \eqref{20b} is $2vnK^2$, and the FLOPS in the dual subgradient calculation in \eqref{eq:20b} is $(2v_mK+1)(n-d)K$. Note that the central operator needs extra $(r-1)(n-d)K$ FLOPS due to \eqref{21c}. However, this is relatively small and negligible compared to the FLOPS of either \eqref{20b} or \eqref{eq:20b}. Therefore, when analyzing the computational cost, we can ignore \eqref{21c} and only consider the computational cost of Group $r_{m}$. Consequently, the cost reduction for the  dual gradient calculation in \eqref{eq:20b} w.r.t. FLOPS are
\begin{subequations}
\begin{align}
    \mathcal{F}_{dt} &= 2vnK^2 - (2v_mK^2+K)(n-d)\nonumber\\
    &= (2vn-2v_m(n-d))K^2 - K(n-d),
    \label{eq:26a}\\
    \mathcal{F}_{dr} &= \frac{(2vn-2v_m(n-d))K^2 - K(n-d)}{2vnK^2}.
    \label{eq:26b}
\end{align}
\end{subequations}

\noindent {\bf{Remark 2:}} For a larger-scale distribution network and a longer time interval $K$ (unnecessarily restricted to valley-filling problem), the dimension reduction $d$ plays a critical role in computational cost reduction. Besides, the dimension reduction makes it possible to decrease the memory cost, therefore taking advantage of the micro controller units embedded in the EV chargers. This will largely facilitate the deployment of the decentralized EV charging control architecture. \hfill $\blacksquare$




\section{Simulation Results}

\subsection{Parameter Configurations}
We consider two scenarios to verify the efficacy and efficiency of the proposed SPMDS-based decentralized EV charging control for valley-filling. First, we introduce some common parameters for both scenarios. The valley-filling period starts from 19:00 to 8:00 next day, with 15-min time intervals and $K=52$ time slots correspondingly. The baseline load is adopted from \cite{liu2017decentralized}. All EVs have the maximum charging power of $6.6$ kW, EVs' random charging requirements vary from $20\%$ to $60\%$ SOC, and the uniform charging efficiency is $\eta_i = 0.9$. The lower bound of the voltage is set to be $\underline{v} = 0.954$, which is slightly higher than the ANSI C84.1 standard to compensate for the LinDistFlow model inaccuracy. The slack node voltage magnitude is $V_0 = 4.16$ kV for both scenarios. 

\subsubsection{Scenario 1: Simulation on the modified IEEE 13-bus Test Feeder} The modified IEEE 13-bus test feeder is adopted from our previous work \cite{liu2017decentralized}.
In this case, 50 EVs are connected at each node except that Nodes 1 and 6 have no load. The grouping and voltage subset selection strategy has been previously shown in TABLE \ref{table_example_13}. The primal step size is set as $\alpha_{i, \ell} = 2.8 \times 10^{-10}$ and the dual step size are set equally as $\beta_{s,\ell} = 1.8$, for $s =1,2,3$. The shrinking parameters are empirically chosen to be $\tau_{\mathcal{U}} =\tau_{\bm{{\lambda}}_1} =\tau_{\bm{\lambda}_2}=\tau_{\bm{\lambda}_3} = 0.98$. 

\begin{figure*}[!htb]
  \centering
  \begin{tabular}[b]{c}
    \includegraphics[width=.45\linewidth,height =.3\linewidth]{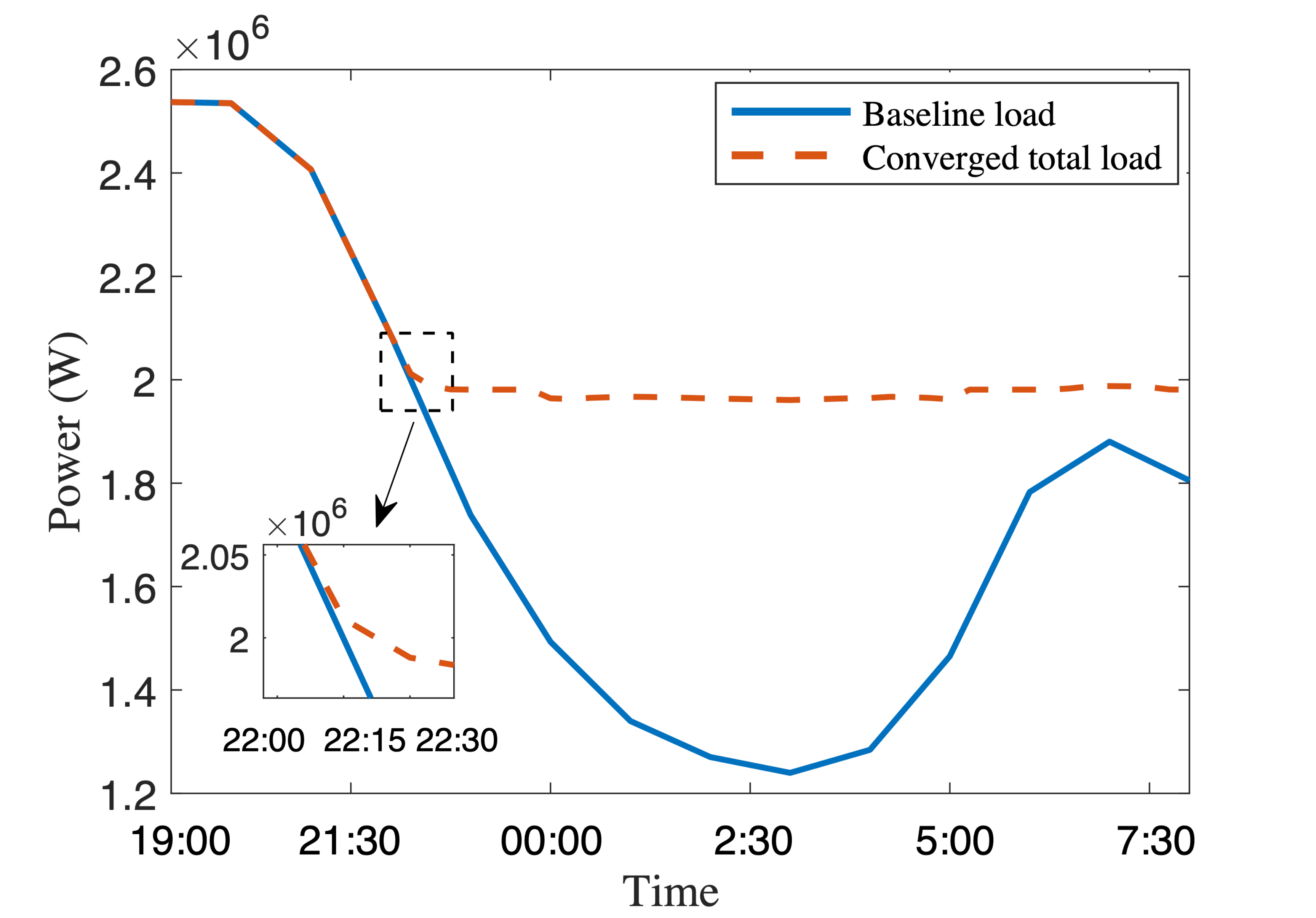} \\
    \small (a)
  \end{tabular} \qquad
  \begin{tabular}[b]{c}
    \includegraphics[width=.45\linewidth,height =.3\linewidth]{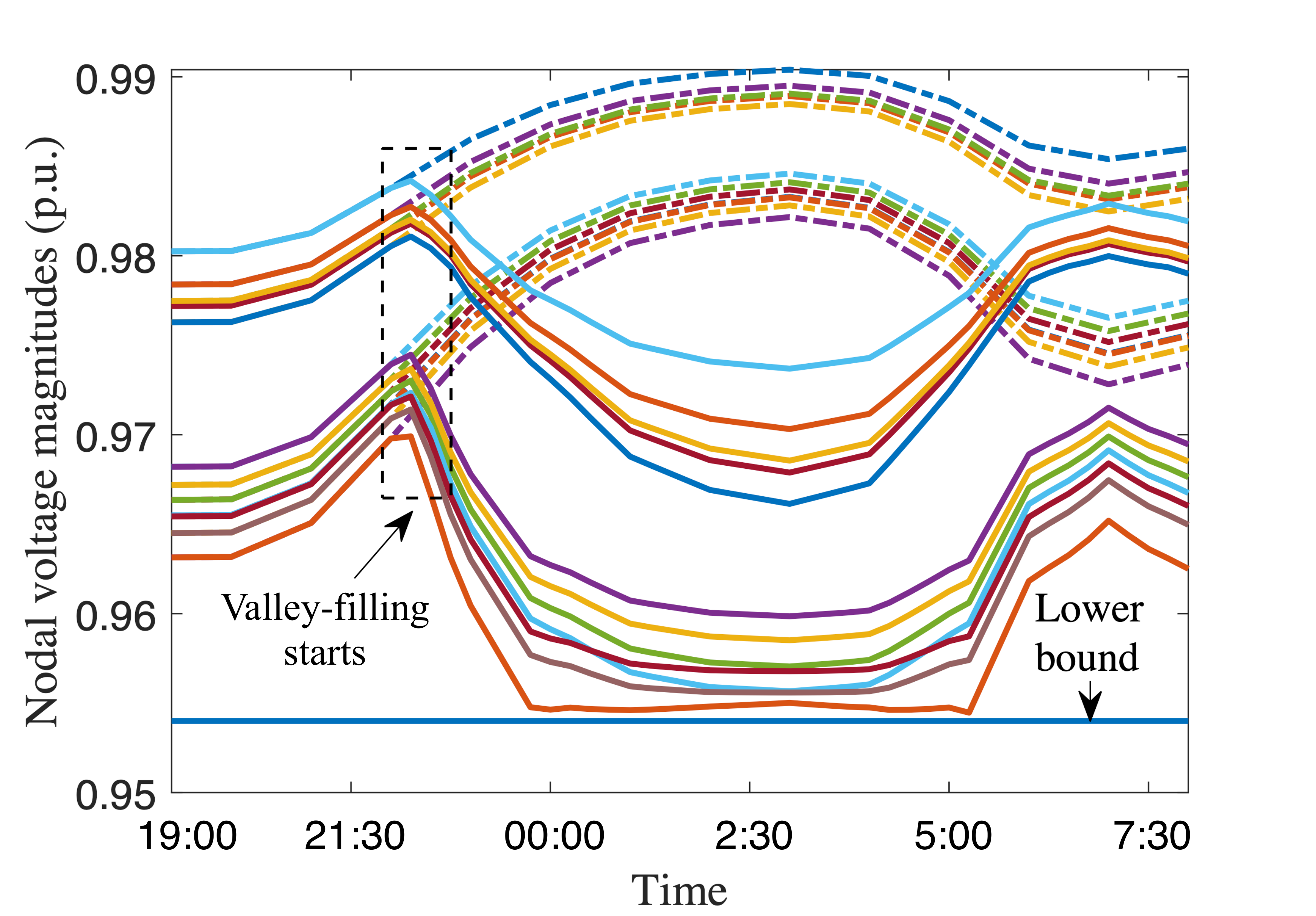} \\
    \small (b)
  \end{tabular}
  \caption{Valley-filling results for 500 EVs on the modified IEEE 13-bus test feeder (20 iterations) (a) Baseline load (solid line) and total load at the 20th iterations (dashed line) (b) Nodal voltage magnitudes of baseline load (dashed lines) and total load (solid lines) }
\label{sim_IEEE_13}
\end{figure*}

The simulation results after 20 iterations are presented in Fig. \ref{sim_IEEE_13}. Fig. \ref{sim_IEEE_13}(a) shows the baseline load in contrast to the total load with the participation of 500 EVs -- the valley caused by the overnight baseline load reduction is filled by the controlled EV charging load between 22:00 and 8:00 the next day. Fig. \ref{sim_IEEE_13}(b) shows that the voltage magnitudes of all nodes are well maintained above the lower bound $0.954$ p.u.

The dimension reduction in this scenario is $d = 5$. As a result, comparing with the SPDS, $\mathcal{F}_{pt} = 27,040$ FLOPS are reduced in the primal subgradient calculation, i.e., $\mathcal{F}_{pr} = 35.76 \%$. For the dual update, Group 3 has the largest number of EVs, i.e., $v_m = 300$. Then the FLOPS reduction in the dual subgradient calculation is $\mathcal{F}_{dt} = 21,090,836$, i.e., $\mathcal{F}_{dr} = 65\%$. 

\subsubsection{Scenario 2: Simulation on the modified IEEE 123-bus Test Feeder}
To verify the scalability of the proposed SPMDS algorithm, we simplified the modified IEEE 123-bus test feeder used in \cite{kersting1991radial,bronzini2011coordination}. In particular, we broke the connection between Node 51 and Node 300 to make it a radial distribution network. In addition, 5 EVs are connected at each node except that Nodes 1 and 6 have no load.

\begin{figure}[!htb]
    \centering
    \includegraphics[width=0.4\textwidth]{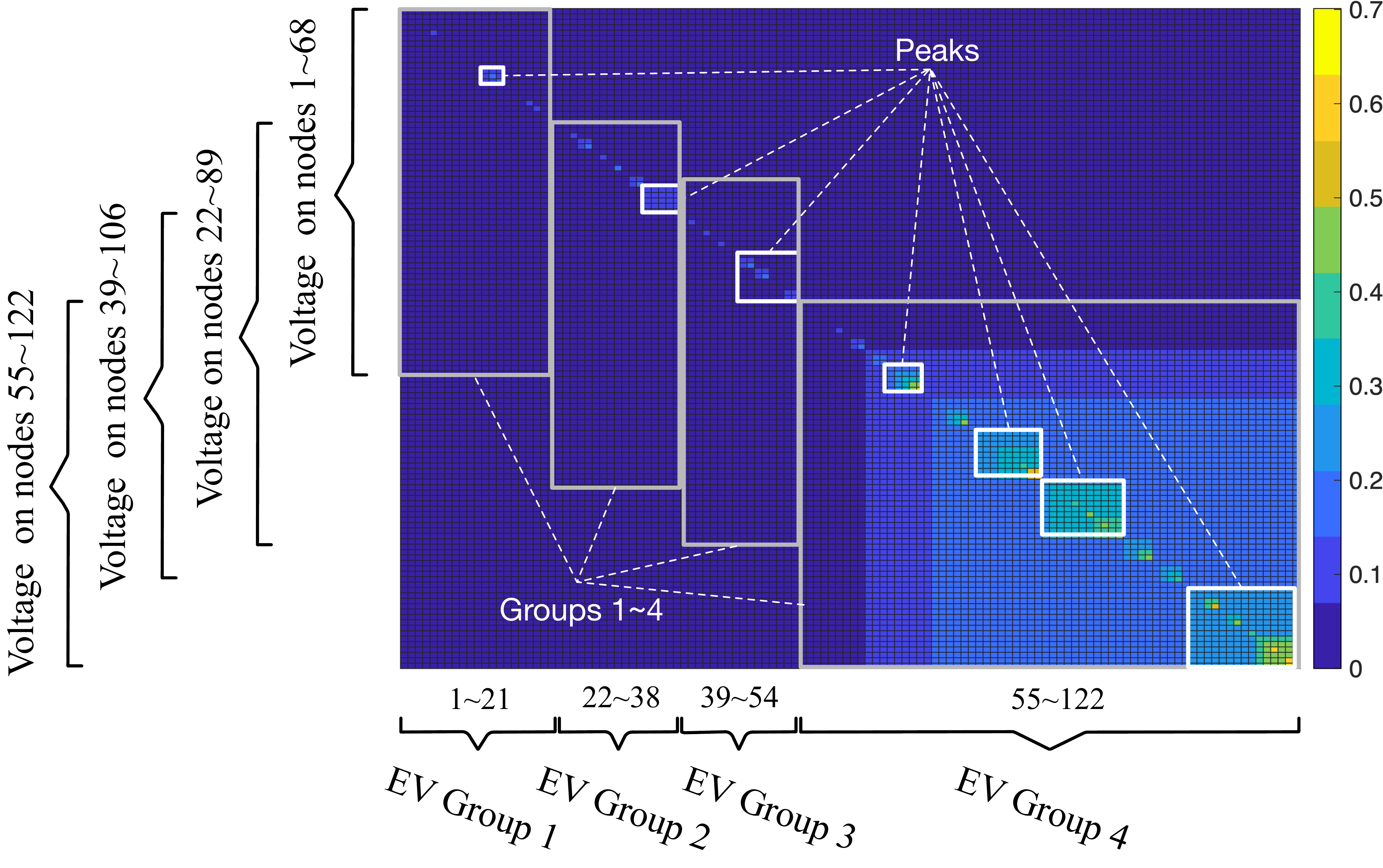}
    \caption{EV grouping and voltage subsets for the modified IEEE 123-bus test feeder (presented via the heatmap of the peak-preserved matrix).}
    \label{sim_IEEE_123_heat}
\end{figure}

The EV grouping and voltage subset selections are presented in Fig. \ref{sim_IEEE_123_heat}  via the heatmap of the peak-preserved matrix where some of the observed ``peaks'' are earmarked by the rectangles. The dimension reduction is set to $d=54$. The total 600 EVs are divided into 4 groups. TABLE \ref{table_example_123} \begin{table}[!htb]
\caption{EV Grouping methodology and voltage subset for each group in the modified IEEE 123-bus test feeder}
\label{table_example_123}
\begin{center}
\begin{tabular}{c|l|l|c}
\hline
Group Name & Location & Voltage subset & Number of EVs\\
\hline
Group 1 & Nodes 1-21    & Nodes  1-68 & 95\\
Group 2 & Nodes 22-38   & Nodes  22-89 & 85\\
Group 3 & Nodes 39-54   & Nodes  39-106& 80\\
Group 4 & Nodes 55-122  & Nodes 55-122 & 340\\
\hline
\end{tabular}
\end{center}
\end{table}presents the details of the grouping and voltage subsets. The primal step size is $\alpha_{i, \ell} = 3 \times 10^{-10}$ and the dual step sizes are $\beta_{s,\ell} = 0.1$, for $s =1,2,3,4$. The shrinking parameters are $\tau_{\mathcal{U}} =\tau_{\bm{{\lambda}}_1}=\tau_{\bm{{\lambda}}_2}=\tau_{\bm{{\lambda}}_3} =\tau_{\bm{{\lambda}}_4} = 0.97$, which is slightly lower than that of the first scenario in exchange for faster convergence. 

\begin{figure*}[!htb]
  \centering
  \begin{tabular}[b]{c}
    \includegraphics[width=.45\linewidth,height =.3\linewidth]{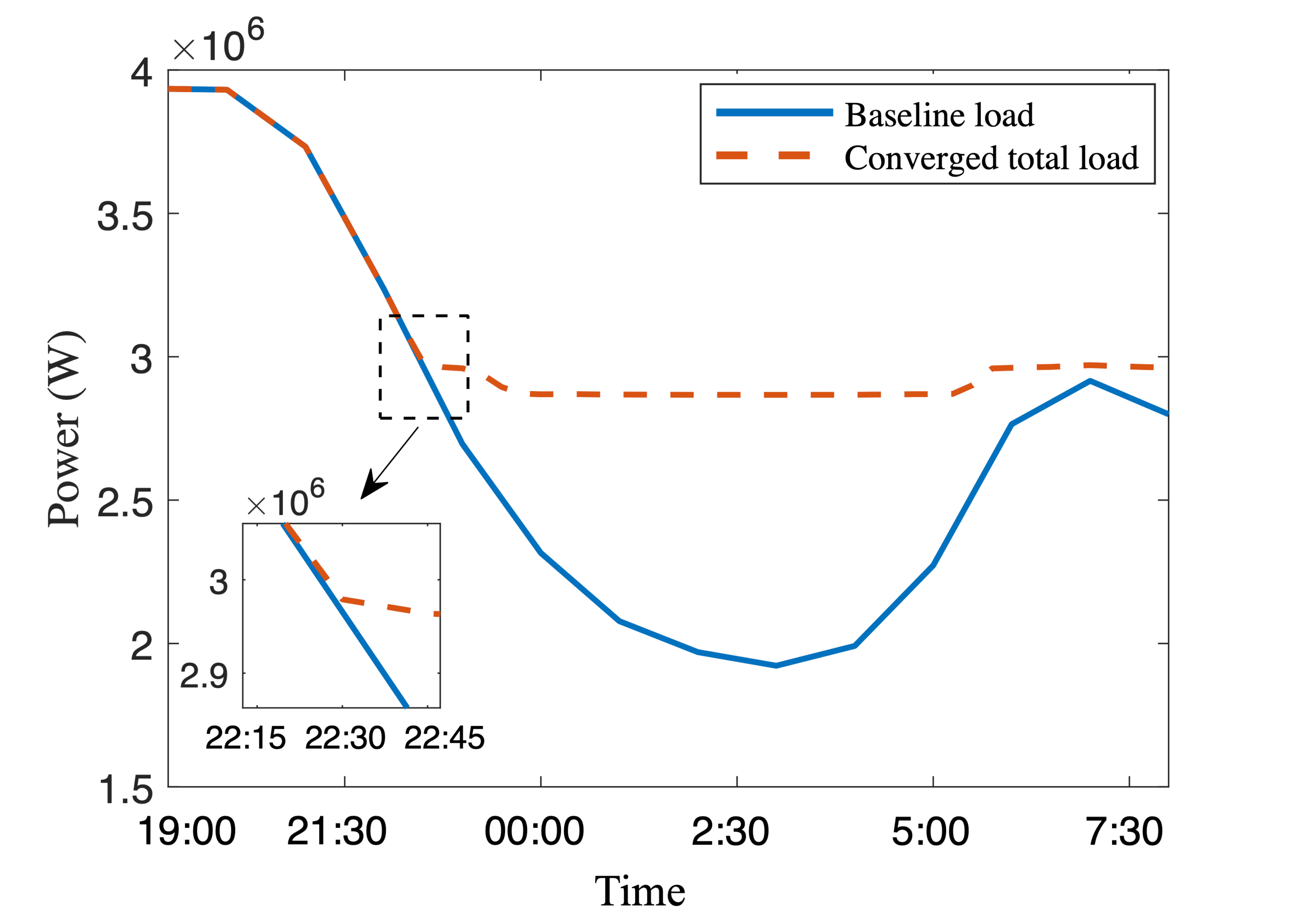} \\
    \small (a)
  \end{tabular} \qquad
  \begin{tabular}[b]{c}
    \includegraphics[width=.45\linewidth,height =.3\linewidth]{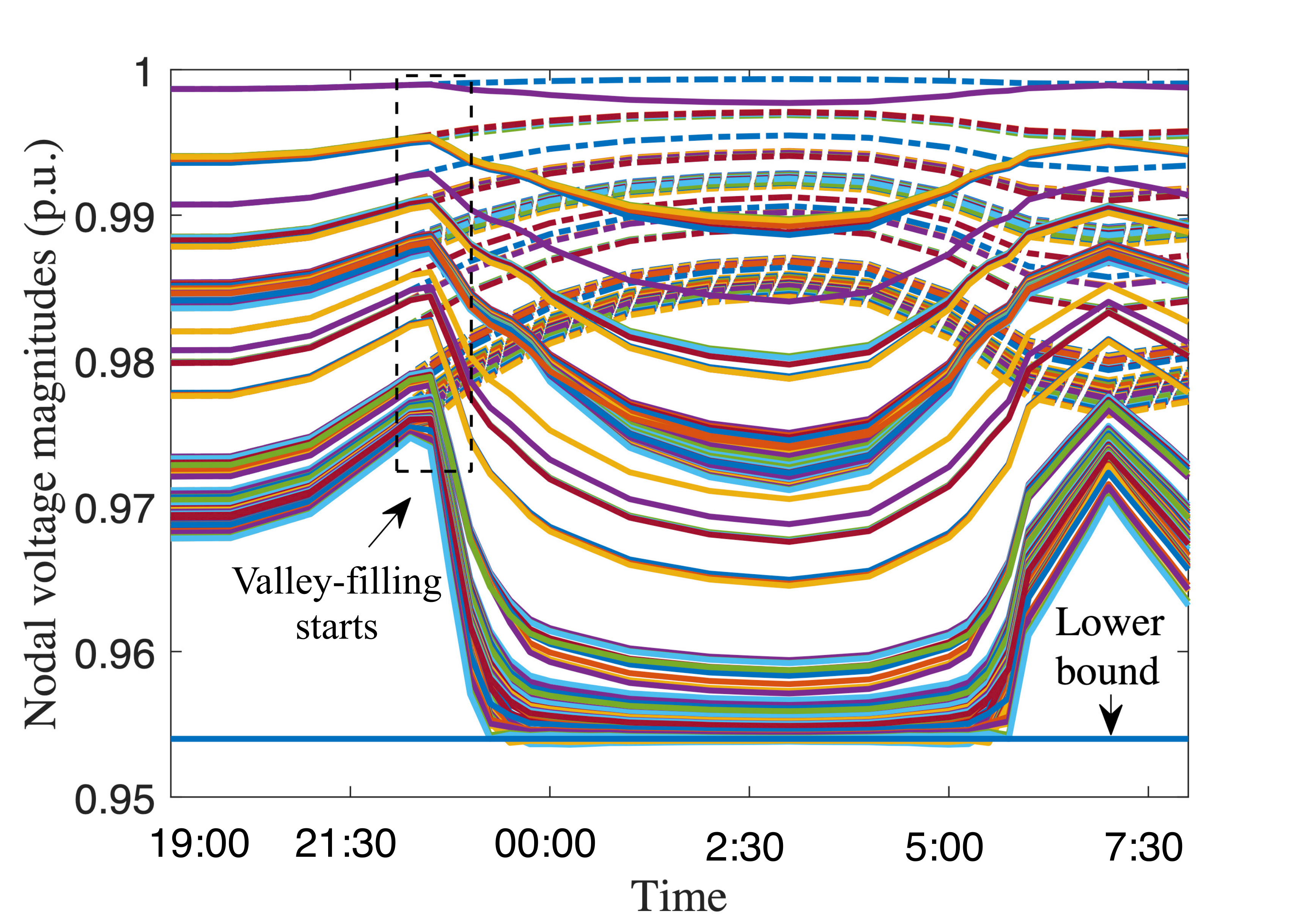} \\
    \small (b)
  \end{tabular}
  \caption{Valley-filling results for 600 EVs on the modified IEEE 123-bus test feeder (40 iterations) (a) Baseline load (solid line) and total load at the 40th iterations (dashed line) (b) Nodal voltage magnitudes of baseline load (dashed lines) and total load (solid lines) }
  \label{sim_IEEE_123}
\end{figure*}

Fig. \ref{sim_IEEE_123}(a) and Fig. \ref{sim_IEEE_123}(b) depict the valley-filling performance and the voltage control of 122 nodes, respectively, which clearly indicate the efficacy of the proposed method. After 40 iterations, though Fig. \ref{sim_IEEE_123}(b) exhibits subtle voltage violations, the converged results are already good enough to satisfy the engineering use. The engineering remedy to the subtle violation could be tightening the constraint bounds, which was already integrated by setting a slightly higher voltage lower bound.

Because of the dimension reduction $d=54$, $\mathcal{F}_{pt} = 292,032$ FLOPS are reduced in the primal subgradient calculation, i.e. $\mathcal{F}_{pr} = 43.56 \%$ FLOPS reduction compared to the full dimensional case. For the dual update, Group 4 has the largest number of EVs, i.e., $v_m = 340$, and $\mathcal{F}_{dt} = 270,829,104$ FLOPS in the dual subgradient calculation are reduced, i.e.,  $\mathcal{F}_{dr} = 68.41\%$ FLOPS reduction compared to the full-dimension case.

\section{Conclusions} 


This paper developed a novel decentralized optimization algorithm, SPMDS, and demonstrated it through establishing a decentralized EV charging control framework for valley-filling. The SPMDS was designed based on a reduced-dimension primal and multi-dual architecture which was achieved by a suite of novel approaches of abstracting and concentrating key topology information of the network (distribution network), grouping agents (EVs), and constructing subsets of the global constraint (voltage) for agent (EV) groups. The SPMDS is scalable w.r.t. the number of agents due to its decentralized architecture, and w.r.t. the network dimension due to its internal primal multi-dual architecture. Simulation results showed the satisfactory valley-filling performance, well regulated network voltage, and the reduced FLOPS in algorithm iterations. SPMDS is not only valid for valley-filling and other optimal power flow problems, it is in fact a powerful scalable tool for generic optimization problems with strongly coupled objective functions and constraints. The novel primal multi-dual architecture also improves the cyber-security of the existing decentralized algorithms. Future work would be investigating the nonuniform network node partitioning, systematic approach of selecting the best dimension reduction $d$, and the convergence analysis. 



\bibliographystyle{IEEEtran}
\bibliography{bibliography}






\end{document}